\documentclass[12pt]{amsart}
\usepackage{amsmath}
\usepackage{amscd}
\usepackage{amsfonts}
\usepackage{amssymb}
\usepackage{amsthm}
\usepackage[all]{xy}

\begin{document}

%New Theorems
\newcommand{\scsltitle}[2]{\par\noindent\textsc{#1}\ \textsl{#2}}
%%%%%%%%%%%%%%%%%%%%%%%%%%%%5
%enumerate style
\renewcommand{\labelenumi}{(\alph{enumi})}

\newtheorem{lemma}{Lemma}[section]
\newtheorem{definition}[lemma]{Definition}
\newtheorem{theorem}[lemma]{Theorem}
\newtheorem{corollary}[lemma]{Corollary}
\newtheorem{fact}[lemma]{Fact}
\newtheorem{proposition}[lemma]{Proposition}
\newtheorem{claim}{Claim}
%%%%%%%%%%%%%%%%%%%%%%
%A,B,C enumerated theorems.
\newtheorem{at}{Theorem}
\renewcommand{\theat}{\Alph{at}}

\newenvironment{example}{\vskip9pt\noindent\addtocounter{lemma}{1}\textbf{Example
\thelemma}}{\hfill $\square$ \par}

\newenvironment{remark}{\vskip9pt\noindent\addtocounter{lemma}{1}\textbf{Remark
\thelemma}}{\hfill $\square$ \par}

%%%%%%%%%%
%\newcounter{counter_enumi}\setcounter{counter_enumi}{0}
\newenvironment{label1}{
\begin{enumerate}
\renewcommand{\labelenumi}{(\arabic{section}.\arabic{lemma}\alph{enumi})}
}{\end{enumerate}\renewcommand{\labelenumi}{(\alph{enumi})}}{}
%%%%%%%%%%%%%%%%%%%%%%%%%%%%%%%%%%%%%%%%%%%%%%%%%
%enumeration of equation by \equation\Lcounter 1a, 1b, 1c etc...
\newcounter{Lcounter}
\renewcommand{\theLcounter}{\theequation\alph{Lcounter}}
\newenvironment{labellist}{\addtocounter{equation}{1}
\begin{list}{{\rm(\arabic{equation}\alph{Lcounter})}}{\usecounter{Lcounter}\leftmargin=5pt}
}{\end{list}}
%%%%%%%%%%%%%%%%%%%%%%%%%%%%%%%%%%%%%%%%%%%%%%%%%%%%%%%%%%
%Eqnarray with enumeration 1a, 1b, 1c, etc...
\newcounter{NumberEquation}

\newenvironment{labeleqnarray}{%
%Initializing NumberEquation
\setcounter{NumberEquation}{\value{equation}}\addtocounter{NumberEquation}{1}
\setcounter{equation}{0}
\renewcommand{\theequation}{\theNumberEquation\alph{equation}}
\begin{eqnarray}%
}%
{\end{eqnarray}\setcounter{equation}{\value{NumberEquation}}
\renewcommand{\theequation}{\arabic{equation}}}
%%%%%%%%%%%%%%%%%%%%%%%%%%%%%%%%%%%%%%%%%%%%%%%%%%%%%%%%%%

%New Commands
\let\origwr=\wre
\def\wre{\mathop{\rm wr}\nolimits}
\renewcommand{\phi}{\varphi}
\newcommand{\lsdp}{\ltimes}
\newcommand{\isom}{\cong}
\newcommand{\invlim}{{\displaystyle \lim_{\longleftarrow}}}
\newcommand{\<}{\left<}
\renewcommand{\>}{\right>}
\newcommand{\Ind}{{\rm Ind}}

\newcommand{\kssigma}{K_s(\mbox{\boldmath $\sigma$})}
\newcommand{\ksrho}{K_s(\mbox{\boldmath $\rho$})}
\newcommand{\gal}[1]{\textnormal{Gal}\left(#1\right)}
\newcommand{\bfsigma}{\mbox{\boldmath$\sigma$}}
\newcommand{\bfrho}{\mbox{\boldmath$\rho$}}
\newcommand{\const}{{\rm const}}
\newcommand{\AwrG}{A\textnormal{wr}G}
\newcommand{\wrG}[1]{#1\textnormal{wr}G}
\newcommand{\awrg}[2]{#1\textnormal{wr}#2}
\newcommand{\nek}{,\ldots,}
\newcommand{\semidirect}{\ltimes}
\newcommand{\eps}{\varepsilon}
\newcommand{\vphi}{\varphi}
\newcommand{\measure}[2]{\mu\left(\left\{#1\mid #2\right\}\right)}
%Braces
\newcommand{\tribra}[1]{\left< #1\right>}
\newcommand{\norbra}[1]{\left( #1\right)}
\newcommand{\limitto}{\mathop{\to}}
%%%%%%%%%%%%%%%%%%%%%%%%%%%%%%%%%%%%%%%%%%%%%%%%%%%%%%%%%%
%from moshe.mac
\def\dotunion{
\def\dotunionD{\bigcup\kern-9pt\cdot\kern5pt}
\def\dotunionT{\bigcup\kern-7.5pt\cdot\kern3.5pt}
\mathop{\mathchoice{\dotunionD}{\dotunionT}{}{}}}

\def\iff{\leftrightarrow}
\def\Iff{\Longleftrightarrow}
\def\inter{\bigcap}
\def\nonempty{\ne\emptyset}
\def\normal{\triangleleft}
\def\notdivide{\nmid}
\def\tensor{\otimes}
\def\union{\bigcup}
%function
\def\ab{{\rm ab}}
\def\and{{\rm and}}
\def\cd{{\rm cd}}
\def\alg{{\rm alg}}
\def\Aut{{\rm Aut}}
\def\gcd{{\rm gcd}}
\def\GL{{\rm GL}}
\def\Hom{{\rm Hom}}
\def\Homeo{{\rm Homeo}}
\def\id{{\rm id}}
\def\Im{{\rm Im}}
\def\Inn{{\rm Inn}}
\def\irr{{\rm irr}}
\def\Ker{{\rm Ker}}
\def\mod{\;\mathop{\rm mod}\hskip0.2em\relax}
\def\ord{{\rm ord}}
\def\Out{{\rm Out}}
\def\PGL{{\rm PGL}}
\def\PSL{{\rm PSL}}
\def\pr{{\rm pr}}
\def\rank{\mathop{\rm rank}\nolimits}
\def\ring{{\rm ring}}
\def\res{{\rm res}}
\def\Res{{\rm Res}}
\def\SL{{\rm SL}}
\def\Spec{{\rm Spec}}
\def\Subg{{\rm Subg}}
%%%%%%%%%%%%%%%%%%%%%%%%%%%%%%%%%%%%%%%%%%%%%%%%%%%%%%%%%%%%%
%FROM moshe.fnt
%cal
\newcommand{\calA}{{\mathcal A}}
\newcommand{\calB}{{\mathcal B}}
\newcommand{\calC}{{\mathcal C}}
\newcommand{\calD}{{\mathcal D}}
\newcommand{\calE}{{\mathcal E}}
\newcommand{\calF}{{\mathcal F}}
\newcommand{\calG}{{\mathcal G}}
\newcommand{\calH}{{\mathcal H}}
\newcommand{\calI}{{\mathcal I}}
\newcommand{\calJ}{{\mathcal J}}
\newcommand{\calK}{{\mathcal K}}
\newcommand{\calL}{{\mathcal L}}
\newcommand{\calM}{{\mathcal M}}
\newcommand{\calN}{{\mathcal N}}
\newcommand{\calO}{{\mathcal O}}
\newcommand{\calP}{{\mathcal P}}
\newcommand{\calQ}{{\mathcal Q}}
\newcommand{\calR}{{\mathcal R}}
\newcommand{\calS}{{\mathcal S}}
\newcommand{\calT}{{\mathcal T}}
\newcommand{\calU}{{\mathcal U}}
\newcommand{\calV}{{\mathcal V}}
\newcommand{\calW}{{\mathcal W}}
\newcommand{\calX}{{\mathcal X}}
\newcommand{\calY}{{\mathcal Y}}
\newcommand{\calZ}{{\mathcal Z}}
%tilde
\newcommand{\agal}{{\tilde a}}       \newcommand{\Agal}{{\tilde A}}
\newcommand{\bgal}{{\tilde b}}       \newcommand{\Bgal}{{\tilde B}}
\newcommand{\cgal}{{\tilde c}}       \newcommand{\Cgal}{{\tilde C}}
\newcommand{\dgal}{{\tilde d}}       \newcommand{\Dgal}{{\tilde D}}
\newcommand{\egal}{{\tilde e}}       \newcommand{\Egal}{{\tilde E}}
\newcommand{\fgal}{{\tilde f}}       \newcommand{\Fgal}{{\tilde F}}
\newcommand{\ggal}{{\tilde g}}       \newcommand{\Ggal}{{\tilde G}}
\newcommand{\hgal}{{\tilde h}}       \newcommand{\Hgal}{{\tilde H}}
\newcommand{\igal}{{\tilde i}}       \newcommand{\Igal}{{\tilde I}}
\newcommand{\jgal}{{\tilde j}}       \newcommand{\Jgal}{{\tilde J}}
\newcommand{\kgal}{{\tilde k}}       \newcommand{\Kgal}{{\tilde K}}
\newcommand{\lgal}{{\tilde l}}       \newcommand{\Lgal}{{\tilde L}}
\newcommand{\mgal}{{\tilde m}}       \newcommand{\Mgal}{{\tilde M}}
\newcommand{\ngal}{{\tilde n}}       \newcommand{\Ngal}{{\tilde N}}
\newcommand{\ogal}{{\tilde o}}       \newcommand{\Ogal}{{\tilde O}}
\newcommand{\pgal}{{\tilde p}}       \newcommand{\Pgal}{{\tilde P}}
\newcommand{\qgal}{{\tilde q}}       \newcommand{\Qgal}{{\tilde Q}}
\newcommand{\rgal}{{\tilde r}}       \newcommand{\Rgal}{{\tilde R}}
\newcommand{\sgal}{{\tilde s}}       \newcommand{\Sgal}{{\tilde S}}
\newcommand{\tgal}{{\tilde t}}       \newcommand{\Tgal}{{\tilde T}}
\newcommand{\ugal}{{\tilde u}}       \newcommand{\Ugal}{{\tilde U}}
\newcommand{\vgal}{{\tilde v}}       \newcommand{\Vgal}{{\tilde V}}
\newcommand{\wgal}{{\tilde w}}       \newcommand{\Wgal}{{\tilde W}}
\newcommand{\xgal}{{\tilde x}}       \newcommand{\Xgal}{{\tilde X}}
\newcommand{\ygal}{{\tilde y}}       \newcommand{\Ygal}{{\tilde Y}}
\newcommand{\zgal}{{\tilde z}}       \newcommand{\Zgal}{{\tilde Z}}
%bar
\newcommand{\gag}[1]{\bar{#1}}
\def\agag{{\bar a}}     \def\Agag{{\gag{A}}}
\def\bgag{{\bar b}}     \def\Bgag{{\gag{B}}}
\def\cgag{{\bar c}}     \def\Cgag{{\gag{C}}}
\def\dgag{{\bar d}}     \def\Dgag{{\gag{D}}}
\def\egag{{\bar e}}     \def\Egag{{\gag{E}}}
\def\fgag{{\bar f}}     \def\Fgag{{\gag{F}}}
\def\ggag{{\bar g}}     \def\Ggag{{\gag{G}}}
\def\hgag{{\bar h}}     \def\Hgag{{\gag{H}}}
\def\igag{{\bar i}}     \def\Igag{{\gag{I}}}
\def\jgag{{\bar \jmath}}    \def\Jgag{{\bar{J}}}
\def\kgag{{\bar k}}     \def\Kgag{{\gag{K}}}
\def\lgag{{\bar l}}     \def\Lgag{{\gag{L}}}
\def\mgag{{\bar m}}     \def\Mgag{{\gag{M}}}
\def\ngag{{\bar n}}     \def\Ngag{{\gag{N}}}
\def\ogag{{\bar o}}     \def\Ogag{{\gag{O}}}
\def\pgag{{\bar p}}     \def\Pgag{{\gag{P}}}
\def\qgag{{\bar q}}     \def\Qgag{{\gag{Q}}}
\def\rgag{{\bar r}}     \def\Rgag{{\gag{R}}}
\def\sgag{{\bar s}}     \def\Sgag{{\gag{S}}}
\def\tgag{{\bar t}}     \def\Tgag{{\gag{T}}}
\def\ugag{{\bar u}}     \def\Ugag{{\gag{U}}}
\def\vgag{{\bar v}}     \def\Vgag{{\gag{V}}}
\def\wgag{{\bar w}}     \def\Wgag{{\gag{W}}}
\def\xgag{{\bar x}}     \def\Xgag{{\gag{X}}}
\def\ygag{{\bar y}}     \def\Ygag{{\gag{Y}}}
\def\zgag{{\bar z}}     \def\Zgag{{\gag{Z}}}

\def\alphagag{{\bar \alpha}}
\def\phigag{{\bar \varphi}}
\def\psigag{{\bar \varphi}}

%Hat
%roman letters with a hat
\def\ahat{{\hat a}}     \def\Ahat{{\hat A}}
\def\bhat{{\hat b}}     \def\Bhat{{\hat B}}
\def\chat{{\hat c}}     \def\Chat{{\hat C}}
\def\dhat{{\hat d}}     \def\Dhat{{\hat D}}
\def\ehat{{\hat e}}     \def\Ehat{{\hat E}}
\def\fhat{{\hat f}}     \def\Fhat{{\hat F}}
\def\ghat{{\hat g}}     \def\Ghat{{\hat G}}
\def\hhat{{\hat h}}     \def\Hhat{{\hat H}}
\def\ihat{{\hat i}}     \def\Ihat{{\hat I}}
\def\jhat{{\hat j}}     \def\Jhat{{\hat J}}
\def\khat{{\hat k}}     \def\Khat{{\hat K}}
\def\lhat{{\hat l}}     \def\Lhat{{\hat L}}
\def\mhat{{\hat m}}     \def\Mhat{{\hat M}}
\def\nhat{{\hat n}}     \def\Nhat{{\hat N}}
\def\ohat{{\hat o}}     \def\Ohat{{\hat O}}
\def\phat{{\hat p}}     \def\Phat{{\hat P}}
\def\qhat{{\hat q}}     \def\Qhat{{\hat Q}}
\def\rhat{{\hat r}}     \def\Rhat{{\hat R}}
\def\shat{{\hat s}}     \def\Shat{{\hat S}}
\def\that{{\hat t}}     \def\That{{\hat T}}
\def\uhat{{\hat u}}     \def\Uhat{{\hat U}}
\def\vhat{{\hat v}}     \def\Vhat{{\hat V}}
\def\what{{\hat w}}     \def\What{{\hat W}}
\def\xhat{{\hat x}}     \def\Xhat{{\hat X}}
\def\yhat{{\hat y}}     \def\Yhat{{\hat Y}}
\def\zhat{{\hat z}}     \def\Zhat{{\hat Z}}

\def\alphahat{{\hat \alpha}}
\def\phihat{{\hat \varphi}}
\def\Phihat{{\hat \Phi}}
\def\Psihat{{\hat \Psi}}

%%%%%%%%%%%%%%%%%%%%%%%%%%%%%%%%%%%%%%%%%%%%%%%%%%%%%%%%
%Mathbb
\newcommand{\bbA}{\mathbb{A}}
\newcommand{\bbB}{\mathbb{B}}
\newcommand{\bbC}{\mathbb{C}}
\newcommand{\bbD}{\mathbb{D}}
\newcommand{\bbE}{\mathbb{E}}
\newcommand{\bbF}{\mathbb{F}}
\newcommand{\bbG}{\mathbb{G}}
\newcommand{\bbH}{\mathbb{H}}
\newcommand{\bbI}{\mathbb{I}}
\newcommand{\bbJ}{\mathbb{J}}
\newcommand{\bbK}{\mathbb{K}}
\newcommand{\bbL}{\mathbb{L}}
\newcommand{\bbM}{\mathbb{M}}
\newcommand{\bbN}{\mathbb{N}}
\newcommand{\bbO}{\mathbb{O}}
\newcommand{\bbP}{\mathbb{P}}
\newcommand{\bbQ}{\mathbb{Q}}
\newcommand{\bbR}{\mathbb{R}}
\newcommand{\bbS}{\mathbb{S}}
\newcommand{\bbT}{\mathbb{T}}
\newcommand{\bbU}{\mathbb{U}}
\newcommand{\bbV}{\mathbb{V}}
\newcommand{\bbW}{\mathbb{W}}
\newcommand{\bbX}{\mathbb{X}}
\newcommand{\bbY}{\mathbb{Y}}
\newcommand{\bbZ}{\mathbb{Z}}
%%%%%%%%%%%%%%%%%%%%%%%%%%%%%%%%%%%%%%%%%%%%%%%%%%%%%%%%%%%%%%
%Math bold
\newcommand{\bfa}{\mbox{\boldmath$a$}}
\newcommand{\bfb}{\mbox{\boldmath$b$}}
\newcommand{\bfc}{\mbox{\boldmath$c$}}
\newcommand{\bfd}{\mbox{\boldmath$d$}}
\newcommand{\bfe}{\mbox{\boldmath$e$}}
\newcommand{\bff}{\mbox{\boldmath$f$}}
\newcommand{\bfg}{\mbox{\boldmath$g$}}
\newcommand{\bfh}{\mbox{\boldmath$h$}}
\newcommand{\bfi}{\mbox{\boldmath$i$}}
\newcommand{\bfj}{\mbox{\boldmath$j$}}
\newcommand{\bfk}{\mbox{\boldmath$k$}}
\newcommand{\bfl}{\mbox{\boldmath$l$}}
\newcommand{\bfm}{\mbox{\boldmath$m$}}
\newcommand{\bfn}{\mbox{\boldmath$n$}}
\newcommand{\bfo}{\mbox{\boldmath$o$}}
\newcommand{\bfp}{\mbox{\boldmath$p$}}
\newcommand{\bfq}{\mbox{\boldmath$q$}}
\newcommand{\bfr}{\mbox{\boldmath$r$}}
\newcommand{\bfs}{\mbox{\boldmath$s$}}
\newcommand{\bft}{\mbox{\boldmath$t$}}
\newcommand{\bfu}{\mbox{\boldmath$u$}}
\newcommand{\bfv}{\mbox{\boldmath$v$}}
\newcommand{\bfw}{\mbox{\boldmath$w$}}
\newcommand{\bfx}{\mbox{\boldmath$x$}}
\newcommand{\bfy}{\mbox{\boldmath$y$}}
\newcommand{\bfz}{\mbox{\boldmath$z$}}

\newcommand{\bfA}{\mbox{\boldmath$A$}}
\newcommand{\bfB}{\mbox{\boldmath$B$}}
\newcommand{\bfC}{\mbox{\boldmath$C$}}
\newcommand{\bfD}{\mbox{\boldmath$D$}}
\newcommand{\bfE}{\mbox{\boldmath$E$}}
\newcommand{\bfF}{\mbox{\boldmath$F$}}
\newcommand{\bfG}{\mbox{\boldmath$G$}}
\newcommand{\bfH}{\mbox{\boldmath$H$}}
\newcommand{\bfI}{\mbox{\boldmath$I$}}
\newcommand{\bfJ}{\mbox{\boldmath$J$}}
\newcommand{\bfK}{\mbox{\boldmath$K$}}
\newcommand{\bfL}{\mbox{\boldmath$L$}}
\newcommand{\bfM}{\mbox{\boldmath$M$}}
\newcommand{\bfN}{\mbox{\boldmath$N$}}
\newcommand{\bfO}{\mbox{\boldmath$O$}}
\newcommand{\bfP}{\mbox{\boldmath$P$}}
\newcommand{\bfQ}{\mbox{\boldmath$Q$}}
\newcommand{\bfR}{\mbox{\boldmath$R$}}
\newcommand{\bfS}{\mbox{\boldmath$S$}}
\newcommand{\bfT}{\mbox{\boldmath$T$}}
\newcommand{\bfU}{\mbox{\boldmath$U$}}
\newcommand{\bfV}{\mbox{\boldmath$V$}}
\newcommand{\bfW}{\mbox{\boldmath$W$}}
\newcommand{\bfX}{\mbox{\boldmath$X$}}
\newcommand{\bfY}{\mbox{\boldmath$Y$}}
\newcommand{\bfZ}{\mbox{\boldmath$Z$}}

\CompileMatrices
\address{
School of Mathematics ~ \\
The Raymond and Beverly Sackler Faculty of Exact Sciences\\
Tel-Aviv University ~ \\
Tel-Aviv, Israel } \email{barylior@post.tau.ac.il}
\title[Diamond Theorem]{Diamond Theorem for a finitely generated Free Profinite Group}
\author{Lior Bary-Soroker}
\date{\today}

\begin{abstract}
We extend Haran's Diamond Theorem to closed subgroups of a
finitely generated free profinite group. This gives an affirmative
answer to Problem 25.4.9 in~\cite{FrJ}.
\end{abstract}

\maketitle

\section*{Introduction}
\noindent Haran's Diamond Theorem for Hilbertian fields roughly
states that every extension $M$ of a Hilbertian field $K$
``captured'' between two Galois extensions $M_1$ and $M_2$ of $K$
is itself Hilbertian \cite{H}, \cite[Thm. 13.8.3]{FrJ}. The
Diamond Theorem has an analog for profinite groups, also due to
Haran:

\begin{at}
Let $m$ be an infinite cardinal. Let $F = \Fhat_m$ be the free
profinite group of rank $m$, \ $M_1$, $M_2$ closed normal
subgroups of $F$, and $M$ a closed subgroup of $F$ satisfying
$M_1\cap M_2\leq M$ and $M_i\not\leq M$ for $i=1,2$. Then $M\cong
\Fhat_m$ \cite[Thm.~25.4.3]{FrJ}.
\end{at}

Problem 25.4.9 of \cite{FrJ} asks for a generalization of Theorem
A to the case where $m$ is finite. A first step toward the
solution of that problem is taken in \cite{Jar}. Proposition 1.3
of \cite{Jar} proves an analog of a theorem of Weissauer for
profinite groups:

\begin{at}
Let $F=\Fhat_e$ with $e\ge2$ an integer, $M$ a closed subgroup of
$F$ of an infinite index, $N$ a closed normal subgroup of $F$
contained in $M$, and $M_0$ an open subgroup of $M$ which does not
contain $N$. Then $M_0 \isom\Fhat_\omega$.
\end{at}

Building on Theorems A and B, we settle here Problem 25.4.9 of
\cite{FrJ} by proving a diamond theorem for free profinite groups
of finite rank:

\begin{at}
Let $F=\Fhat_e$ with $e\ge2$ an integer, $M_1,M_2$ closed normal
subgroups of $F$, and $M$ a closed subgroup of $F$ with
$(F:M)=\infty$, $M_1\cap M_2\le M$, $M_1\not\le M$, and
$M_2\not\le M$. Then $M\isom\Fhat_\omega$.
\end{at}

The proof of Theorem A (at least in the case $m=\aleph_0$) is
reduced to solving a finite embedding problem
$$
(\phi\colon F\to G,\, \alpha\colon A\wre_{G_0} G\to G),
$$
where $G$ is a finite group, \ $A$ is a finite nontrivial group,
$G_0$ is a subgroup of $G$ acting on $A$, and $A\wre_{G_0} G$ is
the twisted wreath product. This embedding problem has a solution
because every finite embedding problem for $\Fhat_\omega$ has a
solution. The same is true in the case $F=\Fhat_e$, with $e$ an
integer, if $e\geq\rank(A\wre_{G_0}G)$. However, in general, this
inequality does not hold.

We observe that $\rank (A\wre_{G_0} G) \leq |G| + \rank (A)$. So,
if we replace $F$ by an open subgroup $E$ containing $M$, then by
Nielsen-Schreier $\rank (E)$ increases (linearly depending on
$(F:E)$). The main problem is that replacing $F$ by $E$ changes
the embedding problem $(\phi,\alpha)$. In this change the order of
$G$ may increase and with it also $|G|+\rank(A)$. The precise
condition when it is possible to choose $E$ such that the rank
condition holds is stated in
Proposition~\ref{conditioned_diamond_theorem}. The condition of
Proposition~\ref{conditioned_diamond_theorem} holds if there are
``many'' subgroups between $F$ and $M$
(Lemma~\ref{composition_lemma}). In this case we say that $M$ is
an ``abundant subgroup of $F$''.

It may happen that there are not enough closed subgroups between
$F$ and one of the subgroups $M$, $MM_1$, or $MM_2$ (in which
case, the corresponding subgroup is called ``sparse'' -
Definition~\ref{def:sparse_abundant}). In this case, the above
proof does not work, so we prove that $M \cong \Fhat_\omega$
directly (Lemma~\ref{Sparse Subgroup}) or by using either Theorem
A or Theorem B.

In the last section we generalize Theorem C to pro-$\calC$ groups,
where $\calC$ is a Melnikov formation of finite groups. We also
transfer Theorem C to the theory of Hilbertian fields and prove
the following result:

\begin{at}
Let $K$ be a PAC field with a finitely generated free absolute
Galois group of rank at least $2$. Let  $M_1$ and $M_2$ be Galois
extensions of $K$. Then every infinite extension $M$ of $K$ in
$M_1M_2$ which is contained neither in $M_1$ nor in $M_2$ is a
Hilbertian field.
\end{at}

\noindent\textsc{Acknowledgement:} The author is grateful to Dan
Haran and Moshe Jarden for their comments and suggestions on
drafts of this manuscript.

%End of File intro.tex

\newpage
\section{Haran's wreath product trick}
\def\st{\mathop{\mid\;}}
\noindent We use twisted wreath product in order to prove the
Diamond Theorem under certain conditions.

\begin{remark}(Twisted Wreath Product).
Let $G$ and $A$ be finite groups and let $G_0$ be a subgroup of
$G$ acting on $A$. We define an action of $G$ on
$$
\Ind_{G_0}^GA =\{f\colon G\to A\st f(\sigma\rho)=f(\sigma)^\rho,\
\sigma\in G,\ \rho\in G_0 \}
$$
by $f^\sigma (\tau) = f(\sigma\tau)$. The corresponding semidirect
product $A\wre_{G_0} G := G\ltimes \Ind_{G_0}^GA$ is called the
twisted wreath product \cite[Def.~13.7.2]{FrJ}. There is a natural
quotient map $A\wre_{G_0} G \to G$ which is defined by
$(\sigma,f)\mapsto \sigma$, $\sigma\in G$, $f\in\Ind_{G_0}^GA$.
For $a\in A$, we define $f_a \in \Ind_{G_0}^G A$ by $f_a(\sigma) =
a^\sigma$ for $\sigma\in G_0$ and $f_a(\sigma) = 1$ for $\sigma\in
G\smallsetminus G_0$. The map $a\mapsto f_a$ embeds $A$ into
$\Ind_{G_0}^G A\leq A\wre_{G_0} G$. Finally, we observe that
$A\wre_{G_0} G$ is generated by $G\cup A$ and thus $\rank
(A\wre_{G_0} G)\leq \rank(G)+\rank(A)\leq |G| +\rank (A)$.
\end{remark}

\begin{proposition}\label{conditioned_diamond_theorem}
Let $e\geq 2$ be an integer. Let $F=\Fhat_e$ be the free profinite
group of rank $e$, $M_1,M_2$ closed normal subgroups of $F$, and
$M$ a closed subgroup of $F$ satisfying $M_1\cap M_2 \leq M$ and
$(F:M)=\infty$.
Suppose that: \\
(*) for every $r\in \bbN$ there exists an open subgroup $E$ of $F$
and an open subgroup $E_0$ of $E$ containing $M$ such that the
following holds:
\begin{labellist}
\item \label{conditioned_diamond_theorem_a}
$(M_i\cap E: M_i\cap E_0) \geq 3$ for $i=1,2$.
\item \label{conditioned_diamond_theorem_b}
$(F:E)\geq r\cdot(E:E_{00})$, where $E_{00} = \bigcap_{\sigma\in
E} E_0^\sigma$.
\end{labellist}
Then $M\cong \Fhat_\omega$.
\end{proposition}

\begin{proof}
We break the proof into five parts.

\scsltitle{Part A:}{Embedding problem for $M$.} As a subgroup of
$\Fhat_e$, $M$ is a projective group of rank at most $\aleph_0$.
Therefore, it suffices to show that every finite split embedding
problem has a solution~\cite[Lemma~24.8.2]{FrJ}. Consider an
embedding problem
\begin{eqnarray}\label{embedding_problem}
(\mu\colon M \to B,\, \beta\colon B\ltimes A \to B)
\end{eqnarray}
in which $A,B$ are finite groups, $B$ acts on $A$, $\beta$ is the
quotient map, and $\mu$ is an epimorphism. We have to find an
epimorphism $\nu\colon M \to B\ltimes A$ with $\beta\circ \nu =
\mu$.

In order to do so, choose an open normal subgroup $D$ of $F$ with
$M\cap D\leq \Ker(\mu)$ and put $r = (F:D) + \rank(A)$. Condition
(*) gives an open subgroup $E$ of $F$ and an open subgroup $E_0$
of $E$ containing $M$ which satisfy
(\ref{conditioned_diamond_theorem_a}) and
(\ref{conditioned_diamond_theorem_b}) with respect to $r$.
Consider the open normal subgroup $L = E_{00} \cap D$ of $E$.
These subgroups satisfy
\begin{eqnarray}\label{eq:c}
(E:L) = (E:E_{00})(E_{00}:L) \leq (E:E_{00})(F:D).
\end{eqnarray}
For $i=1,2$ write $M'_i = M_i\cap E$. Then,
by~(\ref{conditioned_diamond_theorem_a}),
$$
(M'_iML:ML) = (M'_i : M'_i\cap ML)\geq (M_i\cap E : M_i\cap
E_0)\geq 3.
$$
In addition, since $e\geq 2$, it follows
by~(\ref{conditioned_diamond_theorem_b}) and (\ref{eq:c}) that
\begin{eqnarray*}
1 + (F:E) (e-1) &\geq& 1 + (E:E_{00}) \cdot r (e-1)\\
&\geq & (E:E_{00})((F:D)+\rank(A))(e-1)\\
& \geq& (E:L) + \rank(A)
\end{eqnarray*}
Thus,
\begin{labeleqnarray}
\label{GroupConditions} M'_2\not\leq ML \quad \mbox{and} \quad (M'_1 ML:ML)\geq 3,\\
\label{RankCondition}1+(F:E)(e-1) \geq (E:L) + \rank (A).
\end{labeleqnarray}
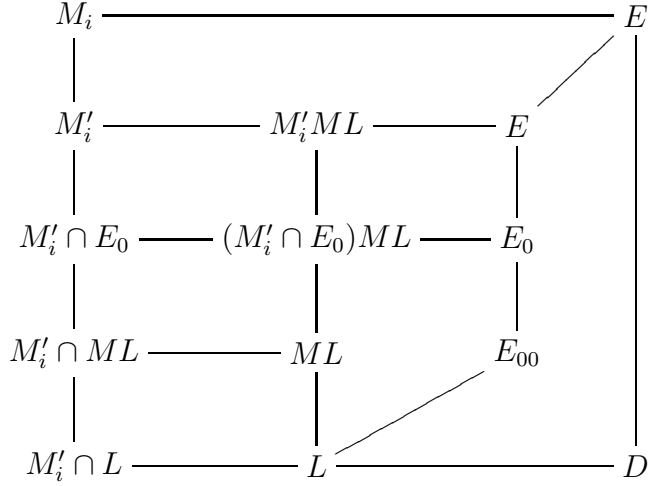
$$
\xymatrix{%
M_i \ar@{-}[rrr]&&&E\\
M'_i \ar@{-}[u] \ar@{-}[r]&M'_iML \ar@{-}[r]& E \ar@{-}[ur] \\
M'_i \cap E_0 \ar@{-}[r]
\ar@{-}[u]&(M'_i \cap E_0)ML \ar@{-}[u] \ar@{-}[r]&E_0 \ar@{-}[u]\\
M'_i \cap ML \ar@{-}[r] \ar@{-}[u]&ML \ar@{-}[u]& E_{00} \ar@{-}[u]\\
M'_i \cap L \ar@{-}[u] \ar@{-}[r]
& L \ar@{-}[u] \ar@{-}[rr] \ar@{-}[ur]&&D \ar@{-}[uuuu]%
}
$$
Since $M\cap L \leq M\cap D \leq \Ker(\mu)$, we can extend $\mu$
to an epimorphism $\phi_1\colon ML\to B$ by $\phi_1(ml) = \mu(m)$
for $m\in M$ and $l\in L$.
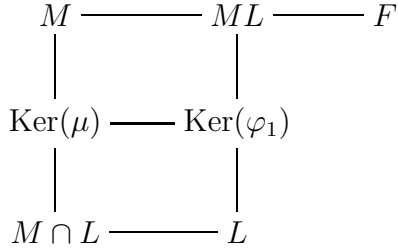
$$
\xymatrix{%
M\ar@{-}[r] &ML \ar@{-}[r]&F\\
\Ker(\mu)\ar@{-}[r]\ar@{-}[u]&\Ker(\phi_1) \ar@{-}[u]\\
M\cap L\ar@{-}[u]\ar@{-}[r]&L\ar@{-}[u] }
$$
Let $G_0=ML/L$ and let $\phi_0\colon ML \to G_0$ be the quotient
map. Then, $G_0=ML/L=\phi_0(ML)=\phi_0(M)$. By definition, $L\leq
\Ker (\phi_1)$, so $\phi_1$ decomposes as in the following
commutative diagram:
\begin{eqnarray*}
\xymatrix{ & ML \ar[d]_{\phi_0} \ar@/^/[dd]^{\phi_1}
\\
G_0\ltimes A \ar[r]^(0.6){\alpha_0} \ar[d]_\rho & G_0
\ar[d]_{\phigag_1}
\\
B\ltimes A \ar[r]^(0.6){\beta} & B }
\end{eqnarray*}
The action of $G_0$ on $A$ is defined via $\phigag_1$. In other
words: $a^\sigma = a^{\phigag_1(\sigma)}$ for $a\in A$, $\sigma\in
G_0$. Also, $\rho|_{G_0} = \phigag_1$ and $\rho|_A = \id_A$.

\scsltitle{Part B:}{The rank and index condition for $E$.} As an
open subgroup of $F$, \ $ML$ is free. If we knew that $\rank
(ML)\geq \rank(B\ltimes A)$, we could find a solution to the
embedding problem $(\phi_1,\beta)$. However, we cannot ensure that
this solution maps $M$ onto $B\ltimes A$. In order to overcome
this difficulty we show in Part C how the embedding problem
$(\phi_0,\alpha_0)$ induces an embedding problem $(\phi\colon F\to
G,\,\alpha\colon A\wre_{G_0} G \to G)$ a solution of which leads
to a solution of (\ref{embedding_problem}). Again, in order to
solve $(\phi,\alpha)$ we need that $\rank(F) \geq \rank
(A\wre_{G_0} G)$. This condition is not necessary fulfilled for
$F$. But it is fulfilled for $E$, as we now show.

Let $G=E/L$ and $G_i = M'_i L/L$ for $i=1,2$. Then $(G:G_0) =
(E:ML)$. By Nielsen-Schreier \cite[Prop. 17.6.2]{FrJ}, $E$ is free
of rank $e' = 1 + (F:E)(e-1)$. Thus, since $\rank(A\wre_{G_0}
G)\leq \rank(G)+\rank(A)$ and by (\ref{RankCondition}),
\begin{eqnarray*}
\rank (A\wre_{G_0}G) &\leq& \rank (G) + \rank(A)\\
&\leq & |G| + \rank(A) \\
&=& (E:L) + \rank (A) \leq  e'
\end{eqnarray*}
Thus, by (\ref{GroupConditions}),
\begin{labeleqnarray}
\label{GroupConditions1}& G_2\not\leq G_0 \quad \mbox{and} \quad
(G_1G_0:G_0)\geq 3.\\
\label{RankCondition1}&\rank(A\wre_{G_0} G) \leq e'.
\end{labeleqnarray}

\scsltitle{Part C:}{Twisted wreath product.} The quotient map
$\phi\colon E\to G$ extends $\phi_0\colon ML \to G_0$. Let
$\alpha\colon A\wre_{G_0} G \to G$ be the quotient map of the
twisted wreath product. By (\ref{RankCondition1}) and by the
freeness of $E$ (\cite[Prop. 17.7.3]{FrJ}), it follows that there
exists an epimorphism $\psi \colon E \to A\wre_{G_0} G$ such that
$\alpha\circ \psi = \phi$.

Every element of $A\wre_{G_0} G = G\ltimes \Ind_{G_0}^G A$ mapped
by $\alpha$ to $G_0$ is in $G_0\ltimes \Ind_{G_0}^GA$. In
particular, $\phi(ML)=G_0$ implies that $\psi(ML)\leq
G_0\ltimes\Ind_{G_0}^GA$.

Define an epimorphism $\pi\colon \Ind_{G_0}^GA \to A$ by $\pi(f) =
f(1)$. Then,
$$
\pi(f)^{\sigma_0} = f(1)^{\sigma_0} = f(\sigma_0) =
f^{\sigma_0}(1) = \pi(f^{\sigma_0}), \quad \forall \sigma_0\in
G_0.
$$
Thus, $\pi$ extends to an epimorphism $\pi\colon G_0\ltimes
\Ind_{G_0}^GA \to G_0\ltimes A$ with $\pi|_{G_0}=\id_{G_0}$. The
following commutative diagram sums up this information.
$$
\xymatrix{%
&&&ML\ar[d]^{\phi_0}\ar[ld]_{\psi|_{ML}}\\
 1\ar[r]&{\Ind_{G_0}^GA}\ar[r]\ar[d]^{\pi}&G_0\ltimes
\Ind_{G_0}^GA \ar[r]^(0.7){\alpha}\ar[d]^{\pi}
& G_0 \ar[r]\ar@{=}[d] & 1\\
1\ar[r]&A\ar[r]&G_0\ltimes A \ar[r]^(0.6){\alpha_0} & G_0 \ar[r] &
1 }
$$
Consider the closed normal subgroup $N = L\cap M'_1\cap M'_2$ of
$E$ and $ML$. By assumption, $M_1\cap M_2\leq M$, so $N\leq M$.

\scsltitle{Part D:}{$\pi(\psi(N)) = A$.} Indeed, $\alpha (\psi(N))
= \phi(N) = 1$, so $\psi(N) \leq \Ind_{G_0}^G A$. In addition the
fact that $N\normal E$ implies that  $\psi(N)\normal A\wre_{G_0}
G$. Hence, $A_1 = \pi(\psi(N))$ is a normal subgroup of
$G_0\ltimes A$ contained in $A$. In particular, $G_0$ preserves
$A_1$.

Assume that $\Agag = A/A_1$ is not trivial. Then $G_0$ acts on
$\Agag$ and we get a commutative diagram
$$
\xymatrix{%
&&&E\ar[d]^{\phi}\ar[ld]_{\psi}\\
1\ar[r]&{\Ind_{G_0}^GA}\ar[r]\ar[d]^{\lambda}& A\wre_{G_0} G
\ar[r]^(0.7){\alpha}\ar[d]^{\lambda} & G \ar[r]\ar@{=}[d] & 1\\
1\ar[r]&\Ind_{G_0}^G\Agag\ar[r]&\Agag\wre_{G_0}G
\ar[r]^(0.6){\alphagag} & G \ar[r] & 1 }
$$
where $\lambda$ is the epimorphism induced by the quotient map
$A\to \Agag$. Now, $\psi(N)\leq \pi^{-1}(A_1) = \{f\in\Ind_{G_0}^G
A\mid f(1)\in A_1\}$ and $\psi(N)$ is preserved by $G$. Thus,
\begin{eqnarray*}
\psi(N) &\le& \inter_{\sigma\in G} \{f\in\Ind_{G_0}^G(A)\mid
f(1)\in A_1\}^\sigma \\
&=&\inter_{\sigma\in G} \{f\in\Ind_{G_0}^G(A)\mid f(\sigma)\in
A_1\} =\Ker(\lambda).
\end{eqnarray*}
Hence,
\begin{eqnarray}
\label{LambdaPsiN1}\lambda(\psi(N))=1.
\end{eqnarray}
For $i=1,2$, put $H_i = \lambda(\psi(M'_i))$. Then $H_i\normal
\Agag\wre_{G_0} G$ and $\alphagag(H_i) = G_i$. By
(\ref{GroupConditions1}) there exists $h_2\in H_2$ with
$\alphagag(h_2) \not\in G_0$. By (\ref{GroupConditions1}) and
\cite[Lemma 13.7.4(a)]{FrJ} there exists $h_1\in H_1$ such that
$\alphagag(h_1) = 1$ and $[h_1,h_2]\neq 1$. For $i=1,2$ choose
$x_i\in M'_i$  such that $\lambda (\psi(x_i)) = h_i$. Then
$\phi(x_1) = \alphagag(h_1) =1$ and thus $x_1\in L$. Hence,
$$
[x_1,x_2]\in [L,M'_2]\cap [M'_1,M'_2]\leq L\cap M'_1\cap M'_2 = N.
$$
It follows from (\ref{LambdaPsiN1}) that
$$
[h_1,h_2] =[\lambda(\psi(x_1)),\lambda(\psi(x_2))]
=\lambda(\psi[x_1,x_2])\in\lambda(\psi(N))=1,
$$
in contradiction to the choice of $h_1$. Consequently,
$\pi(\psi(N)) = A$.

\scsltitle{Part E:}{A solution of the original embedding problem.}
The maps we have defined so far give the following commutative
diagram:
$$
\xymatrix{%
&&&ML\ar[d]_{\phi_0}\ar[ld]_{\psi|_{ML}}\ar@/^0.7pc/[ddd]^{\phi_1} \\
 1\ar[r]&{\Ind_{G_0}^GA}\ar[r]\ar[d]_{\pi}&G_0\ltimes
\Ind_{G_0}^GA \ar[r]^(0.7){\alpha}\ar[d]_{\pi}
& G_0 \ar[r]\ar@{=}[d] & 1\\
1\ar[r]&A\ar[r]\ar@{=}[d] &G_0\ltimes A\ar[d]_{\rho}
\ar[r]^(0.6){\alpha_0}
& G_0 \ar[d]_{\phigag_1}\ar[r] & 1\\
1\ar[r]& A \ar[r] &B\ltimes A \ar[r]^(0.6)\beta &B\ar[r]&1
 }
$$
In particular, $\beta\circ\rho\circ\pi\circ\psi|_M =
\phi_1|_M=\mu$. By Part D,
$$
\rho(\pi(\psi(M)))\geq \rho(\pi(\psi(N)))=\rho(A)=A.
$$
Hence, $\rho(\pi(\psi(M)))=B\ltimes A$. Consequently,
$\rho\circ\pi\circ\psi|_M$ is a solution of the embedding
problem~(\ref{embedding_problem}), as desired.
\end{proof}

The next lemma replaces Condition (*) by a more convenient one. We
use the next observation in the lemma and in the rest of the
paper: Let $E_0$ be an open subgroup of a profinite group $F$.
Then $E_{00}=\bigcap_{\sigma\in F}E_0^\sigma$ is an open normal
subgroup of $F$ and $(F:E_{00})\leq (F:E_0)!$. Indeed, the action
of $F$ on the right cosets of $E_0$ in $F$ induces a homomorphism
of $F$ into the group of permutations $S_n$, where $n=(F:E_0)$,
whose kernel is exactly $E_{00}$.

\begin{lemma}\label{Preperation_Lemma}
Let $F = \Fhat_e$ with $e\geq 2$ an integer, $M_1,M_2$ closed
normal subgroups of $F$, and $M$ a closed subgroup of $F$
satisfying $M_1\cap M_2 \leq M$ and $(M_i:M_i\cap M) = \infty$ for
$i=1,2$. Suppose $F$ has for each $s\in\bbN$ open subgroups $E_1
\leq E$ containing $M$ such that $(F:E)\geq s \cdot (E:E_1)!$ and
for each $i\in\{1,2\}$ either
\begin{labellist}
\item $M_i\leq E$ or \label{Preperation_Lemma_a}
\item $M_iE_1 = F$ and $(E:E_1)\geq 3$.\label{Preperation_Lemma_b}
\end{labellist}
Then $M\cong \Fhat_\omega$.
\end{lemma}

\begin{proof}
We prove that Condition (*) in
Proposition~\ref{conditioned_diamond_theorem} is fulfilled.

Let $r\in \bbN$. Since $(M_i:M_i\cap M)=\infty$, $F$ has an open
subgroup $H$ containing $M$ such that $(M_i : M_i\cap H) \geq 3$
for $i=1,2$. Put $s=r\cdot(F:H)!$. The assumption of the lemma
gives open subgroups $E_1\leq E$ containing $M$ such that
$(F:E)\geq s \cdot (E:E_1)!$ and for each $i\in\{1,2\}$ either
(\ref{Preperation_Lemma_a}) or (\ref{Preperation_Lemma_b}) holds.
Set $E_0 = H\cap E_1$, $E_{00} = \bigcap_{\sigma\in E}
E_0^\sigma$, $E_{11} = \bigcap_{\sigma\in E} E_1^\sigma$ and
$H_{00} = \bigcap_{\sigma\in F}H^\sigma$. Then $H_{00}\cap
E_{11}\leq \bigcap_{\sigma\in E}(H^\sigma\cap E_1^\sigma) =
\bigcap_{\sigma\in E} E_0^\sigma = E_{00}$.
$$
\xymatrix{%
H_{00}\ar@{-}[r]\ar@{-}[dd] & H\ar@{-}[r]\ar@{-}[d] &
F\ar@{-}[d]|{\mbox{$E$}} \\
& E_0\ar@{-}[r]\ar@{-}[dl]|{\mbox{$E_{00}$}}&E_1\ar@{-}[d]\\
H_{00}\cap E_{11}\ar@{-}[rr]&&E_{11}
 }
$$
Hence
\begin{eqnarray*}
(E:E_{00}) & \leq & (E:H_{00}\cap E_{11})\\
& = & (E:E_{11})(E_{11} : H_{00} \cap E_{11})\leq
(E:E_{11})(F:H_{00})\\
& \leq & (E:E_1)!(F:H)!\leq  \frac{1}{s} (F:E)(F:H)! =
\frac1r(F:E)
\end{eqnarray*}
This proves (\ref{conditioned_diamond_theorem_b}).

In order to prove (\ref{conditioned_diamond_theorem_a}) we first
assume that $M_i\leq E$. Then, since $E_0\leq H$,
$$
(M_i\cap E: M_i\cap E_0) \geq (M_i:M_i\cap H)\geq 3.
$$

Now assume that $M_i E_1 = F$ and $(E:E_1)\geq 3$. Then $(M_i \cap
E) E_1 = E$, so
$$
(M_i\cap E:M_i\cap E_0)\geq (M_i\cap E:M_i\cap E_1) = (E:E_1) \geq
3.
$$
It follows from Proposition~\ref{conditioned_diamond_theorem} that
$M\cong \Fhat_\omega$.
\end{proof}

\newpage
\section{Sparse and abundant subgroups}
We can prove that the condition of Lemma~\ref{Preperation_Lemma}
is satisfied only if there are ``many'' subgroups between $F$ and
$M$. For example, if $MM_1$ is ``abundant'' in $F$ in the sense of
the following definition. Luckily, if $MM_1$ is ``sparse'' in $F$
then it is isomorphic to $\Fhat_\omega$ (Lemma~\ref{Sparse
Subgroup}). Use of Haran's Diamond Theorem for profinite groups of
infinite rank yields in this case that $M$ itself is isomorphic to
$\Fhat_\omega$(Theorem~\ref{Diamond_Theorem}).

\begin{definition}
[Sparse and Abundant subgroups]\label{def:sparse_abundant} A
closed subgroup $M$ of a profinite group is called {\bf sparse} if
for all $m,n\in\bbN$ there exists an open subgroup $K$ of $F$
containing $M$ such that $(F:K) \geq m$, and for every proper open
subgroup $L$ of $K$ containing $M$ we have $(K:L)\geq n$. In
particular, $(F:M)=\infty$.

A closed subgroup $M$ of $F$ with $(F:M)=\infty$ which is not
sparse is said to be {\bf abundant}.
\end{definition}

\noindent\textsl{Example of a sparse subgroup:} Let $F = \prod
\bbZ/p\bbZ$ where $p$ runs over all prime numbers. Then $1$ is
sparse in $F$.

\begin{lemma}\label{Improvment Lemma}
Let $M$ be a sparse subgroup of a profinite group $F$. Then for
each open subgroup $H$ of $F$ containing $M$ and for all $m,n\in
\bbN$ there exists an open subgroup $K$ of $H$ containing $M$ with
$(F:K)\geq m$ and $(K:L)\geq n$ for each open proper subgroup $L$
of $K$ containing $M$.
\end{lemma}

\begin{proof}
Let $H$ be an open subgroup of $F$ containing $M$, and let $m,n\in
\bbN$. Since $M$ is sparse in $F$, $F$ has an open subgroup $K$
containing $M$ such that $(F:K)\geq m$ and $(K:L) \geq 1 +
\max\{n,(F:H)\}$ for each proper open subgroup $L$ of $K$
containing $M$, in particular $(K:L)>(F:H)$. Then, since $(K:K\cap
H)\leq(F:H)$ and $M\leq K\cap H$, it follows that $K\cap H = K$,
i.e.~$K\leq H$.
\end{proof}

\begin{corollary}\label{Sparse_unde_Taking_Subgroups}
Let $M\leq H\leq F$ be profinite groups with $H$ open in $F$.
Then, $M$ is sparse in $F$ if and only if $M$ is sparse in $H$.
\end{corollary}

\begin{lemma}\label{Sparse Subgroup}
Let $F = \Fhat_n$ with $2\leq n \leq \aleph_0$ and $M$ a subgroup
of $F$. If $M$ is sparse in $F$, then $M\cong \Fhat_\omega$.
\end{lemma}

\begin{proof}
The rank of $M$ is at most $\aleph_0$, so, by Iwasawa, it suffices
to prove that every finite embedding problem $(\phi\colon M \to A,
\alpha \colon B\to A)$ for $M$ has a solution \cite[Cor.
24.8.3]{FrJ}. Indeed, choose an open subgroup $D\normal F$ with
$D\cap M\leq \Ker(\phi)$ and set $H = MD$. Then $H$ is open in $F$
and $\phi$ extends to an epimorphism $\phi'\colon H \to A$ by
$\phi'(md)=\phi'(m)$ for $m\in M$ and $d\in D$. Since $M$ is
sparse, there is an open subgroup $K$ of $H$ containing $M$ which
has no proper open subgroup $L$ containing $M$ and satisfying
$(K:L)\leq |B|$ (Lemma~\ref{Improvment Lemma}). Moreover, if $n$
is finite, we can choose $K$ such that $1+(F:K)(n-1)\geq\rank(B)$.
The latter inequality is obvious when $n=\aleph_0$.

By Nielsen-Schreier, $K$ is a free profinite group of rank at
least $\rank(B)$. Hence, there exists an epimorphism $\gamma\colon
K\to B$ with $\alpha\circ \gamma = \phi'|_K$ (\cite[Prop.
17.7.3]{FrJ} if $n$ is finite and \cite[Thm. 24.8.1]{FrJ}
otherwise). By the choice of $K$ in the preceding paragraph,
$\Ker(\gamma) M = K$. Hence,
$\gamma(M)=\gamma(\Ker(\gamma)M)=\gamma(K)=B$. Consequently,
$\gamma|_M$ is a solution of the embedding problem
$(\phi,\alpha)$.
\end{proof}

\begin{lemma}\label{Abundant_Subgrp}
Let $M$ be an abundant subgroup of a profinite group $F$. Then for
each $s\in \bbN$ there exist open subgroups $E_1\leq E$ of $F$
containing $M$ such that $(F:E)\geq s\cdot (E:E_1)!$ and
$(E:E_1)\geq 3$.
\end{lemma}

\begin{proof}
By definition, there exist $m,n\in\bbN$ such that for every open
subgroup $K$ of $F$ containing $M$ with $(F:K)\geq m$ there exists
an open subgroup $L$ containing $M$ such that $1<(K:L)\leq n$.

Let $s\in \bbN$. Since $(F:M)=\infty$, there exists an open
subgroup $K$ of $F$ containing $M$ with $(F:K)\geq \max\{s\cdot
n!, s\cdot 4!,m\}$. By assumption, $K$ has an open subgroup $L$
containing $M$ such that $1<(K:L)\leq n$.

If $(K:L)\geq 3$ the subgroups $E=K$ and $E_1 = L$ satisfy the
conclusion of the lemma. Otherwise, $(K:L) = 2$. By assumption $L$
has an open subgroup $L_0$ containing $M$ such that $1<(L:L_0)\leq
n$. If $(L:L_0)\geq 3$, the subgroups $E=L$ and $E_1 = L_0$
satisfy the conclusion of the lemma. Otherwise, $(L:L_0) = 2$ and
$(K:L_0) = 4$, so $E=K$, $E_1=L_0$ satisfy the conclusion of the
lemma.
\end{proof}

\begin{lemma}\label{composition_lemma}
Let $F = \Fhat_e$ with $e\geq 2$, $M_1,M_2$ open normal subgroups
of $F$, and $M$ a closed subgroup satisfying $M_1\cap M_2 \leq M$
and $(M_i:M_i\cap M)=\infty$ for $i=1,2$. In addition, assume that
at least one of the following conditions holds:
\begin{labellist}
\item \label{composition_lemma_a}
$(F:MM_1M_2) = \infty$.
\item \label{composition_lemma_b}
$(F:MM_1M_2)<\infty$ and $MM_1$ is abundant in $F$.
\item \label{composition_lemma_c}
$(F:MM_1M_2)<\infty$ and $MM_2$ is abundant in $F$.
\item \label{composition_lemma_d}
$(F:(MM_1)\cap (MM_2))<\infty$ and $M$ is abundant in $F$.
\end{labellist}
Then $M \cong \Fhat_\omega$.
\end{lemma}

\begin{proof}
Let $s\in \bbN$. By Lemma~\ref{Preperation_Lemma} it suffices to
find open subgroups $E_1\leq E$ of $F$ containing $M$ such that
$(F:E)\geq s \cdot (E:E_1)!$ and for each $i\in\{1,2\}$ either
\begin{labellist}
\item
\label{cond:1} $M_i\leq E$ or
\item \label{cond:2}$M_i E_1
= F$ and $(E:E_1)\geq 3$.
\end{labellist}
We distinguish between the four cases:

\scsltitle{Case A:}{$(F:MM_1M_2)=\infty$.} Choose an open subgroup
$E$ of $F$ containing $MM_1M_2$ such that $(F:E)\geq s$. Put $E_1
= E$. Then $M_1,M_2\leq E$ and $(F:E)\geq s\cdot (E:E_1)!$.

\scsltitle{Case B:}{$(F:MM_1M_2)<\infty$ and $MM_1$ is abundant in
$F$.} By Corollary~\ref{Sparse_unde_Taking_Subgroups}, we can
replace $F$ by $MM_1M_2$ in order to assume that $F = MM_1M_2$; it
suffices to proof $(9)$ for this $F$. By
Lemma~\ref{Abundant_Subgrp}, $F$ has open subgroups $E_1 < E$
containing $MM_1$ such that $(F:E)\geq s\cdot (E:E_1)!$ and
$(E:E_1)\geq 3$. Then, $M_1\leq E$ and $E_1M_2 = F$.

\scsltitle{Case C:}{$(F:MM_1M_2)<\infty$ and $MM_2$ is abundant in
$F$.} Exchange the indices $1$ and $2$. Then the result follows
from Case B.

\scsltitle{Case D:}{$(F: (MM_1)\cap (MM_2))<\infty$ and $M$ is
abundant in $F$.} Let $F' = (MM_1)\cap (MM_2)$. By
Corollary~\ref{Sparse_unde_Taking_Subgroups} $M$ is abundant in
$F'$.
%
%, in particular $(F':M)=\infty$.  For $i=1,2$,
%$$
%(M_i: M_i\cap M) = (MM_i:M) \geq (F':M) = \infty.
%$$
Put $M'_1 = M_1\cap F'$ and $M'_2 = M_2\cap F'$. Then $MM'_1 =
MM'_2 = F'$. By assumption $(MM_i:M)=(M_i:M_i\cap M)=\infty$.
Therefore, since $(M_i:M'_i) < \infty$ it follows that
$(M'_i:M'_i\cap M) = \infty$.
$$
\xymatrix{%
M_1\ar@{-}[r]&M_1M \ar@{-}[r]&F\\
M'_1\ar@{-}[r]\ar@{-}[u]& F'\ar@{-}[r]\ar@{-}[u]& M_2M\ar@{-}[u]\\
&M'_2\ar@{-}[r]\ar@{-}[u]&M_2\ar@{-}[u] }
$$
Replace $F$ by $F'$, $M_1$ by $M'_1$, and $M_2$ by $M'_2$, if
necessary, to assume that $MM_1 = F$ and $MM_2 = F$; again it
suffices to prove $(9)$ for this $F$. Lemma~\ref{Abundant_Subgrp}
gives open subgroups $E_1\leq E$ of $F$ containing $M$ with
$(F:E)\geq s\cdot (E:E_1)!$ and $(E:E_1)\geq 3$. Those subgroups
satisfy (\ref{cond:2}), i.e., $M_1E_1 = F$ and $M_2 E_1 = F$.
\end{proof}

\begin{theorem}[Diamond Theorem]\label{Diamond_Theorem}
Let $F=\Fhat_e$ with $e\ge2$ an integer, $M_1,M_2$ closed normal
subgroups of $F$, and $M$ a closed subgroup of $F$ with
$(F:M)=\infty$, $M_1\cap M_2\le M$, $M_1\not\le M$, and
$M_2\not\le M$. Then $M\isom\Fhat_\omega$.
\end{theorem}

\begin{proof}
If $(M_2:M_2\cap M)<\infty$, then $(MM_2:M)=(M_2:M\cap
M_2)<\infty$, so $M$ is a proper open subgroup of $MM_2$ with
$M_2\not\leq M$. By Theorem~B, $M\cong \Fhat_\omega$. The same
argument gives that if $(M_1:M_1\cap M) < \infty$, then $M\cong
\Fhat_\omega$. Thus, we may assume that
\begin{eqnarray}\label{Both_Infinity}
(M_i:M_i\cap M)=\infty, \text{ for } i=1,2.
\end{eqnarray}

If $(F:MM_1M_2) =\infty$, then (\ref{composition_lemma_a}) holds.
Hence, by Lemma~\ref{composition_lemma}, $M\cong \Fhat_\omega$.
Thus, we may assume that $(F:MM_1M_2)<\infty$. By
Nielsen-Schreier, $MM_1M_2$ is a finite ranked profinite group of
rank at least $2$. Thus, we can replace $F$ by $MM_1M_2$ to assume
that $F= MM_1M_2$.

Let $F' = MM_1\cap MM_2$. We distinguish between two case.

\scsltitle{Case A:}{$(F:F')<\infty$.} If $M$ is abundant in $F$,
then (\ref{composition_lemma_d}) holds. Hence, by
Lemma~\ref{composition_lemma}, $M\cong \Fhat_\omega$. If $M$ is
sparse in $F$, then by Lemma~\ref{Sparse Subgroup},
$M\cong\Fhat_\omega$.

\scsltitle{Case B:}{$(F:F') = \infty$.} Then either $(F:MM_1) =
\infty$ or $(F:MM_2)=\infty$. Without loss we may assume that
$(F:MM_1) = \infty$.

If $MM_1$ is abundant in $F$, then (\ref{composition_lemma_b})
holds. Hence, by Lemma~\ref{composition_lemma}, $M\cong
\Fhat_\omega$.

Assume therefore that $MM_1$ is sparse in $F$. By
Lemma~\ref{Sparse Subgroup}, $MM_1 \cong \Fhat_\omega$. Put $M'_2
= MM_1\cap M_2$. Then $M_1, M'_2\normal MM_1$, $M_1\cap M'_2\leq
M$ and $M_1\not \leq M$. If $M'_2\not\leq M$, so by Theorem A,
$M\cong \Fhat_\omega$. Otherwise, $M'_2\leq M$ and thus
$M=MM_1\cap MM_2 = F'$. In particular
$(F:MM_2)=(MM_1:M)=(M_1:M_1\cap M)=\infty$ (by
(\ref{Both_Infinity})).
$$
\xymatrix{MM_1\ar@{-}[r]&F\\
M=F'\ar@{-}[u]\ar@{-}[r]&MM_2\ar@{-}[u]\\
M'_2\ar@{-}[u]\ar@{-}[r] &M_2\ar@{-}[u] }
$$
If $MM_2$ is abundant in $F$, then (\ref{composition_lemma_c})
holds. Hence, by Lemma~\ref{composition_lemma}, $M\cong
\Fhat_\omega$. If $MM_2$ is sparse in $F$, then $M$ is sparse in
$MM_1$ (because $MM_1/M_2' \cong F/M_2$). Hence, by Lemma
\ref{Sparse Subgroup}, $M\cong \Fhat_\omega$.

In each case we have proven that $M\cong\Fhat_\omega$.
\end{proof}

\newpage
\section{Applications}
\noindent As in the infinite rank case \cite[Thm 25.4.4]{FrJ}, it
is possible to generalize the Diamond Theorem to the category of
pro-$\calC$ groups for each Melnikov formation $\calC$. This
generalization (Theorem 3.1 below) settles Problem~25.4.9 in
\cite{FrJ}.

Recall that a family $\calC$ of finite groups is called a
\textbf{Melnikov formation} if it is closed under taking
quotients, normal subgroups, and extensions. We write
$\Fhat_m(\calC)$ for the free pro-$\calC$ group of rank $m$.

\begin{theorem}\label{Melnikov}
Let $\calC$ be a Melnikov formation, $2\leq e <\infty$, and
$M_1,M_2,M$ be closed subgroups of $F = \Fhat_e(\calC)$. Suppose
that $M$ is a pro-$\calC$ group, $(F:M)=\infty$, $M_1,M_2\normal
F$, $M_1\cap M_2 \leq M$ but $M_1\not\leq M$ and $M_2\not\leq M$.
Then $M\cong \Fhat_{\omega}(\calC)$
\end{theorem}

\begin{proof}
Let $\Fhat = \Fhat_e$ and $\Nhat = M_{\Fhat}(\calC)$ the
intersection of all open normal subgroups $K$ of $\Fhat$ with
$\Fhat/K\in \calC$. Lemma 17.4.10 of \cite{FrJ} gives an
epimorphism $\phi\colon \Fhat \to F$ with $\Ker(\phi) = \Nhat$.
Put $\Mhat_1 = \phi^{-1} (M_1)$, $\Mhat_2 = \phi^{-1} (M_2)$ and
$\Mhat = \phi^{-1} (M)$. Then $\Mhat_1,\Mhat_2\normal \Fhat$,
$\Nhat\leq \Mhat_1\cap \Mhat_2\leq \Mhat$ but $\Mhat_1\not \leq
\Mhat$ and $\Mhat_2\not\leq \Mhat$. By
Theorem~\ref{Diamond_Theorem}, $\Mhat \cong \Fhat_\omega$.

In order to prove that $M \cong \Fhat_\omega(\calC)$ it suffices
now to show that $\Nhat=M_{\Mhat}(\calC)$ \cite[Lemma
17.4.10]{FrJ}. Indeed, let $L$ be an open normal subgroup of
$\Mhat$ with $\Mhat/L\in \calC$. Then $\Nhat / L\cap \Nhat \cong
L\Nhat /L$ and $L\Nhat/L\normal \Mhat/L$, so $\Nhat/L\cap \Nhat\in
\calC$. By \cite[Lemma 17.4.10]{FrJ}, $L\cap \Nhat = \Nhat$, hence
$\Nhat\leq L$. It follows that $\Nhat \leq M_{\Mhat}(\calC)$. On
the other hand, $\Mhat/\Nhat\cong M$ is a pro-$\calC$ group, so
$M_{\Mhat}(\calC)\leq \Nhat$. Consequently,
$\Nhat=M_{\Mhat}(\calC)$, as claimed.
\end{proof}

There is a notable special case of Theorem~\ref{Melnikov}

\begin{corollary}
Let $\calC$ be a Melnikov formation of finite groups, $2\leq e
<\infty$, and $M_1,M_2$ closed normal subgroups of
$\Fhat_e(\calC)$. Suppose none of the groups $M_1$ and $M_2$ is
contained in the other and at least one of them has an infinite
index. Then $M_1\cap M_2 \cong \Fhat_\omega(\calC)$.
\end{corollary}

Haran's Diamond Theorem has been originally proved for Hilbertian
fields: Let $K$ be a Hilbertian field, $M$ a separable extension
of $K$, and $M_1,M_2$ Galois extensions of $K$. If $M\not\subseteq
M_1$ and $M\not\subseteq M_2$ but $M\subseteq M_1M_2$, then $M$ is
Hilbertian \cite[Thm.~13.8.3]{FrJ}. In particular this theorem
holds for $\omega$-free PAC fields, because they are Hilbertian
\cite[Cor. 27.3.3]{FrJ}. On the other hand, if $K$ is just PAC and
$e$-free for $2\leq e<\infty$, then $K$ is not Hilbertian (because
the rank of the absolute Galois group of a Hilbertian field has
infinite index). Nevertheless, the theorem is true also for PAC
$e$-free fields:

\begin{theorem}
Let $K$ be a PAC $e$-free field with $2\leq e<\infty$ and
$M,M_1,M_2$ separable extensions of $K$. Suppose that
$[M:K]=\infty$, $M \not\subseteq M_1$, $M \not\subseteq M_2$, but
$M\subseteq M_1M_2$. Then $M$ is Hilbertian.
\end{theorem}

\begin{proof}
Use Theorem~\ref{Diamond_Theorem} for $\gal{K}$ and its closed
subgroups $\gal{M}$, $\gal{M_1}$, $\gal{M_2}$ to get that
$\gal{M}\cong \Fhat_\omega$. By Ax-Roquette, $M$ is PAC \cite[Cor
11.2.5]{FrJ}. Hence, by Roquette, $K$ is Hilbertian
\cite[Cor.~27.3.3]{FrJ}.
\end{proof}

\newpage

\end{document}